\documentstyle{amsppt}
\NoBlackBoxes

\TagsOnRight

\def\cal{\Cal}

\def\EE{{\cal E}}

\def\FF{{\cal F}}

\def\Z{{\Bbb Z}}
\def\C{{\Bbb C}}
\def\R{{\Bbb R}}

\def\n{\noindent}
\def\del{{\partial}}
\def\delbar{{\overline\partial}}

\rightheadtext{coisotropic submanifolds} \leftheadtext{
Yong-Geun Oh }

\topmatter
\title
Geometry of coisotropic submanifolds in symplectic manifolds
and K\"ahler manifolds
\endtitle
\author
Yong-Geun Oh\footnote{Partially supported by the NSF Grant \#
DMS-0203593, the grant of 2000 Korean Young Scientist Prize
and Vilas Research Award of University of Wisconsin
\hskip10.5cm\hfill}
\endauthor
\address
Department of Mathematics, University of Wisconsin, Madison, WI
53706 ~USA \& Korea Institute for Advanced Study, Seoul, Korea
\endaddress

\abstract The first purpose of this paper is to generalize the
well-known Maslov indices of maps of open Riemann surfaces with
boundary lying on Lagrangian submanifolds to maps with boundary
lying on coisotropic submanifolds in symplectic manifolds. For
this purpose, we first define the notion of {\it Maslov loops} of
coisotropic Grassmanians and their indices. Then we introduce the
notions of {\it transverse Maslov bundle} of coisotropic
submanifolds, and {\it gradable coisotropic submanifolds}. We then
define {\it graded coisotropic submanifolds} and the {\it
coisotropic Maslov index} of the maps with boundary lying on such
graded coisotropic submanifolds, which reduces to the standard
Maslov index of disc maps for the case of Lagrangian submanifolds.
The second purpose is to study the geometry of coisotropic
submanifolds in K\"ahler manifolds. We introduce the notion of the
leafwise mean curvature form and transverse canonical bundle of
coisotropic submanifolds and study various geometric properties
thereof. Finally we combine all these to define the notion of
{\it special coisotropic submanifolds} for the case of Calabi-Yau manifolds,
and prove various consequences on their properties of the coisotropic
Maslov indices.

\endabstract

\endtopmatter

\bigskip

\document

\head \bf \S 1. Introduction \endhead

A triple $(X,\omega,J)$ is called an almost K\"ahler manifold where
$(X,\omega)$ is symplectic and $J$ is compatible to $\omega$ in that
linear form defined by
$$
g(X,Y): = \omega(X, JY)
$$
is positive definite. The {\it mean curvature one form} $\alpha_L$
of the Lagrangian submanifold $L \subset (X,g)$ is then defined by
$$
\alpha_L: = (\vec H \rfloor \omega)|_{TL} \tag 1.1
$$
Morvan [Mo] and Dazord [Da] proved that the mean curvature one
form, of the Lagrangian submanifold $L \subset (X,\omega, J)$
becomes closed when $(X,\omega, J)$ is Einstein-K\"ahler, i.e.,
$J$ is integrable and its Ricci form $\rho$ satisfies
$$
\rho = \lambda\omega, \quad \text{for } \, \lambda\in \R. \tag 1.2
$$
More precisely, Dazord [Da] proved the identity
$$
d\alpha_L = i^*\rho \tag 1.3
$$
on general K\"ahler manifolds from which together with (1.2)
follows closedness of $\alpha_L$ in the Einstein-K\"ahler case. We
will call the {\it mean curvature class} the real one dimensional
cohomology class $[\alpha_L] \in H^1(L;\R)$ of the Lagrangian
submanifold $L \subset (X,\omega,J)$. In fact, Morvan (in the case
of $\C^n$) and Dazord (in the Calabi-Yau, i.e., in the Ricci-flat
case) proved that the  mean curvature one form
$\frac{1}{\pi}\alpha_L$ is the one form whose de Rham cohomology
class is an integral class which represents the well-known
Maslov class [Ar] of Lagrangian submanifolds in symplectic
geometry. The latter measures the rotation index of the {\it
angle} of the tangent plane of the Lagrangian submanifolds which
illustrates an interesting interplay between symplectic and
Riemannian geometry of Lagrangian submanifolds in
Einstein-K\"ahler manifolds. In this respect, the author [Oh2]
also proved that the mean curvature class $\alpha_L$ is invariant
under the Hamiltonian isotopy of $(X,\omega)$ (see [Oh2] for its
precise meaning).

One purpose of the present paper is to generalize all of the above
facts on Lagrangian submanifolds on K\"ahler manifolds to {\it
coisotropic submanifolds} $Y \subset (X,\omega)$: A submanifold
$Y$ is called coisotropic if
$$
(TY)^\omega \subset TY \tag 1.4
$$
where $(TY)^\omega$ is the {\it symplectic orthogonal complement}
to $TY$ defined by
$$
(T_yY)^\omega = \{v \in T_yX \mid  \omega(v, \cdot) \equiv 0 \}.
\tag 1.5
$$
Coisotropic submanifolds have played some important role in
symplectic geometry in relation to generalizing calculus of
Lagrangian submanifolds to some corresponding calculus in Poisson
manifolds in relation to the geometric quantization (see e.g.,
[We]). Recently they attracted some physicist's attention [KaOr]
in an attempt to {\it correct and complete} Kontsevich's
homological mirror symmetry proposal in which inclusion of
coisotropic $D$-branes plays an important role. For the case of
K\"ahler manifolds, we introduce the notion of {\it leafwise mean
curvature vector} of a coisotropic submanifold $Y \subset
(X,\omega, J)$ and attempt to prove all the analogs in this
coisotropic case to the above mentioned relationships between the
Maslov index and the mean curvature class that are valid for
Lagrangian submanifolds. It turns out that in a suitable {\it
foliated} or {\it leafwise} context, all the coisotropic analogs
to the above mentioned geometric properties of Lagrangian
submanifolds can be proved.

We also introduce the notion of {\it Maslov loops} of coisotropic
Grassmanians
$$
\Gamma_k(\R^{2n},\omega_0) = \{C \in Gr_{n+k}(\R^{2n}) \mid
C^{\omega_0} \subset C, \, \dim\ker\omega_0|_C = n-k \}
\tag 1.6
$$
and their indices. Note that when $k = 0$, $C$ become Lagrangian
and $k = n$, it becomes the full space $\R^{2n}$. This class is
{\it not} an homotopy invariant unlike the Lagrangian case but
enjoys certain symplectic invariance property (see section 2 and
3). Then we introduce the notions of {\it gradable coisotropic
submanifolds} and {\it graded coisotropic submanifolds} in general
symplectic manifolds, for which we also define an index for any
(disc) map $w: (D,\del D) \to (X,Y)$.

In section 7, we will also introduce the notion of {\it
special coisotropic submanifolds} in the Calabi-Yau case which
generalize the notion of special Lagrangian submanifolds. And we
will derive various consequence on the coisotropic Maslov index of
disc maps with boundary lying on special coisotropic submanifolds.

Finally in section 8, we analyze the symplectic $\Pi$-transverse
curvature $F_\Pi$ introduced in [OP1] with respect to the orthogonal
splitting $\Pi: TY = T\FF \oplus N_J\FF$
in the case of K\"ahler manifolds $(X,\omega,J)$ and relate the
curvature with the classical Levi form for the case of hypersurfaces
$Y \subset (X,\omega,J)$.

The present work is a by-product of the joint works [OP1,2] with
Jae-Suk Park. We would like to thank Jae-Suk Park for exciting
collaboration on the coisotropic $D$-branes. We also thank
Wei-Dong Ruan for his interest on this work and for pointing out a
couple of imprecise points in our calculations in the precious
version of this paper.

\head \bf \S 2. Maslov loops of coisotropic Grassmanians and
their indices
\endhead

In this sub-section, we introduce the notion of Maslov loops of
coisotropic subspaces and their associated indices.

Denote the set of coisotropic subspace of rank $2k$ in
the symplectic vector space $(\R^{2n},\omega_0)$ by
$$
\Gamma_k(\R^{2n}, \omega_0) = :\{ C \in Gr_{n+k}(\R^{2n}) \mid
C^\omega \subset C, \, \dim C^\omega = n-k \}.
$$
From the definition, for any coisotropic subspace we have the
canonical flag,
$$
0 \subset C^\omega \subset C \subset \R^{2n}.
$$
Combining this with the standard complex structure $j$ on $\R^{2n}
\cong \C^n$, we have the splitting
$$
C = H_C \oplus C^\omega
$$
where $H_C$ is the complex subspace of $C$. The following
proposition is not difficult check, whose proof we omit and leave
to [OP1].

\proclaim{Proposition 2.1} Let $ 0 \leq k \leq n$ be fixed. The
unitary group $U(n)$ acts transitively on $\Gamma_k$. The
corresponding homogeneous space is given by
$$
\Gamma_k(\R^{2n}, \omega_0) \cong U(n)/U(k) \times O(n-k)
$$
where $U(k) \times O(n-k) \subset U(n)$ is the isotropy group of
the coisotropic subspace $\C^k \oplus \R^{n-k} \subset \C^n$. In
particular we have
$$
\dim \Gamma_k(\R^{2n}, \omega_0) = \frac{(n+3k +1)(n-k)}{2}.
$$
\endproclaim
Now we have the following commutative diagram of short exact
sequences
$$
\matrix  0 & \to H & \to C & \to C/C^\omega \to 0 \\
\quad & \quad \downarrow &\quad
\downarrow &\downarrow &\quad \\
0 & \to H & \to \C^n & \to \C^n/H \to 0
\endmatrix
\tag 2.1
$$
Suppose $\gamma: S^1 \to \Gamma_k(\R^{2n},\omega_0)$
be a continuous loop and consider the flag
$$
0 \subset \gamma(\theta)^\omega \subset \gamma(\theta) \subset
\R^{2n}.
$$
Considering
$$
S_\gamma: = C_\gamma^\omega \oplus jC_\gamma^\omega , \quad
C_\gamma^\omega(t) : = \gamma(\theta(t))^\omega
$$
we associate to the each coisotropic loop $\gamma$ the pair
$(S_\gamma, L_\gamma)$ of loops of symplectic vector space
$S_\gamma(\theta)$ of dimension $2(n-k)$ and its Lagrangian
subspace $L_\gamma(\theta)=C_\gamma^\omega$ for $\theta \in S^1$.
We denote by $(C^\omega)^\perp \subset C^*$ the set of
annihilators of $C^\omega$ in the dual space $C^*$ of $C$.

We will need the following general discussion on the coisotropic
subbundle of symplectic vector bundles.  Let $(E,\sigma,J) \to N$
be a symplectic vector bundle with a compatible complex structure
$J$ on the fiber. Let $C \subset E$ be a coisotropic subbundle and
$C^\sigma$ the kernel bundle of $C$. Denote by $C_\sigma =
C/C^\sigma$ the quotient symplectic vector bundle and $J_\sigma$
the induced complex structure. We then identify $C_\sigma$ with
the orthogonal complement of $C^\sigma$ in $C$. We have the
splitting
$$
C_\sigma \otimes_{J_\sigma} \C =C^{1,0}_{\sigma,J} \oplus C^{0,1}_{\sigma,J}.
$$
with respect the complex structure $J$. We recall that as a
complex vector bundle $(C_\sigma, J_\sigma)$ and $C^{0,1}_\sigma
= C^{1,0}_{\sigma,J}$
are canonically isomorphic.  We suppress the $J_\sigma$
from notations and write the splitting
$$
C_\sigma^* \otimes \C = \Lambda^{1,0}(C^*_\sigma) + \Lambda^{0,1}(C^*_\sigma)
$$
associated to $C_\sigma^*$. Note that $C_\sigma^*$ is naturally
isomorphic to $(C^\sigma)^\perp \subset C^*$, the set of annihilators
of $C^\sigma$.

\definition{Definition 2.2} \roster
\item We denote the determinant line
bundle of $\Lambda^{1,0}(C_\sigma^*)$ by
$$
K_C=K_{(C;E,\Omega,J)}: =
\Lambda^{k,0}(C_\sigma^*) = \Lambda^k(\Lambda^{1,0}(C_\sigma^*))
$$
and call it the {\it transverse canonical bundle} of $C$.
We also denote by $\Omega^{k,0}(C_\sigma^*)$ the set of smooth sections
thereof.  We call any section $\zeta \in \Omega^{k,0}(C_\sigma^*)$ of unit norm
a {\it transverse complex volume form} of $C$.
\item We call the square bundle $K_C^{\otimes 2}$ the {\it transverse
Maslov line bundle} of $C \subset (E,\Omega,J)$.
\endroster
\enddefinition

We note that since the real line bundle $(\det_\R(C^\sigma))^{\otimes 2}$
is trivial and $\det_\R E$ is always trivial.
In general, the transverse Maslov bundle $K_C^{\otimes 2}$ is not
trivial as a complex vector bundle.

\definition{Definition 2.3} Let $C \subset (E,\Omega)$ be a
coisotropic subbundle of nullity $n-k$ or of rank $2k$. \roster
\item
We call the coisotropic subbundle $C \subset (E,\sigma)$ {\it
gradable} if $K_C^{\otimes 2}$ is trivial, for a (and so any)
compatible almost complex structure.
\item
For a given compatible almost complex structure $J$ on
$(E,\Omega)$, we call a section $\zeta \in \Gamma(K_C^{\otimes
2})$ a {\it transverse Maslov section} of $C$. We call the pair
$(C,\sigma)$ a {\it graded} coisotropic subbundle of $E$.
\endroster
\enddefinition

\definition{Examples 2.4}
\roster
\item Any Lagrangian subbundle is gradable in this sense since in that case the
transverse Maslov bundle is just a scalar as the transverse space is trivial.
\item Any coisotropic subbundle over the base $N = S^1$ is
gradable since any complex line bundle over $S^1$ is trivial.
\item Consider $(X,\omega,J)$  is Calabi-Yau, or more generally
any symplectic manifold $(X,\omega)$ with $2c_1(X,\omega) = 0$.
Then for any coisotropic submanifold $Y \subset (X,\omega)$,
the coisotropic subbundle $TY \subset (TX|_Y,\omega)$ is gradable
(see section 7).
\endroster
\enddefinition

Now we go back to the case $E=S^1 \times (\R^{2n},\omega_0)$ the
trivial symplectic vector bundle over $S^1$ with the standard
complex structure $j$, with a map $\gamma: S^1 \to
\Gamma_k(\R^{2n},\omega_0)$ which corresponds to a coisotropic
subbundle of $E$. This is a special case where $N = S^1$. In this
case, the real bundle $(C^\sigma)^{\otimes 2} \to S^1$ is always
trivial. Since any symplectic vector bundle is trivial over $S^1$,
both $C_\sigma^*$ and $(\det_\C C_\sigma^*)^{\otimes 2}$ are also
trivial ($\C^n$ is Calabi-Yau!).

\definition{Definition 2.5}
\roster
\item
We denote
$$
K_\gamma = \Lambda^{k,0}(S_\gamma)
$$
and call the {\it transverse canonical bundle} of $\gamma$ (with
respect to the complex structure induced from the standard complex
structure of $\C^n$).  We call any section $\zeta$ of
$K_\gamma^{\otimes 2}$ a {\it transverse Maslov section}
of $\gamma$.
\item We call the pair $(\gamma;\zeta)$ a {\it Maslov loop}.
\endroster
\enddefinition

Now we will associate an integer to each given Maslov loop $(\gamma;\zeta)$.

We fix a global section $(\Omega^k_\R)^{\otimes 2}(C_\gamma^\omega)$ which
we write as
$$
(e_{k+1} \wedge \cdots \wedge
e_n)^{\otimes 2} \tag 2.2
$$
for any local orthonormal frame
$$
\{e_{k+1}, \cdots, e_n\}  \tag 2.3
$$
of $C_\gamma^\omega$. It is obvious that (2.2) is independent of
the choice of orthonormal frame. However $e_{k+1}\wedge \cdots
\wedge e_n$ is not globally defined in general. In case
$C_\gamma^\omega$ is oriented, one can make this globally defined
considering only the oriented frames. We denote by $\{e_1, \cdots,
e_{k+1}, \cdots, e_n, f_1, \cdots, f_{k+1}, \cdots, f_n\}$ any
Darboux frame of $\C^n$, that extends the frame $\{e_{k+1},\cdots,
e_n\}$ and respects the orthogonal splitting
$$
\gamma(\theta) = S_\gamma(\theta) \oplus C_\gamma^\omega(\theta).
$$
More precisely, $\{e_1, \cdots, e_k, f_1, \cdots, f_k\}$ defines
a Darboux orthonormal frame of $S_\gamma$
and $\{e_{k+1}, \cdots, e_n, f_{k+1}, \cdots, f_n\}$
one for $C_\gamma^\omega \oplus jC_\gamma^\omega$.
Denote $L_\gamma = C_\gamma^\omega \subset S_\gamma$ which is a
Lagrangian subbundle of $S_\gamma$.

Noting that $S_\gamma \cong L_\gamma \otimes \C \cong
(S_\gamma)^{1,0}$, we choose the unitary frame
$\{u_1, \cdots, u_k\}$ of $(S_\gamma)^{1,0}$ given as
$$
u_j = e_j + if_j
$$
associated to the Darboux frame $\{e_1, \cdots, e_k, f_1,
\cdots, f_k \}$. We denote by $\{\theta^1, \cdots, \theta^n \}$
the unitary frame dual to $\{u_1, \cdots, u_n\}$.

We know that the square of
$$
u_{k+1} \wedge \cdots \wedge u_n
$$
defines a global section of $K_\gamma^{\otimes 2}$.
Pairing this with the standard complex volume form
$$
dz = dz_1 \wedge \cdots \wedge dz_n,
$$
it gives rise to a natural transverse Maslov section of $\gamma$
$$
\zeta_{(\gamma;dz)}: =
\Big((u_{k+1} \wedge \cdots \wedge u_n) \rfloor dz)\Big)^{\otimes 2}. 
\tag 2.4
$$
The following is easy to check whose proof we omit.

\proclaim{Lemma 2.6} (2.4) does not depend  on the choice of the
orthonormal frame $\{e_{k+1}, \cdots, e_n\}$ of $C_\gamma^\omega$
or its extended orthonormal Darboux frame $\{e_1, \cdots, e_n,
f_1, \cdots, f_n\}$, but depends only on $\gamma$.
\endproclaim

Now let $\zeta$ be any given transverse Maslov section of $\gamma$.
Then we can write
$$
\zeta(\theta) = g(\theta) \zeta_{(\gamma;dz)}(\theta)
\tag 2.5
$$
for a well-defined function $g = g_{(\gamma;\zeta)}: S^1 \to S^1$.

\definition{Definition 2.7} We define the index $\mu(\gamma;\zeta)$
of the Maslov loop $(\gamma;\zeta)$ by
$$
\mu(\gamma;\zeta): = \text{deg} (g_{(\gamma;\zeta)}). \tag 2.6
$$
\enddefinition

\definition{Remark 2.8} \roster
\item
Note that in the Lagrangian case there is no transverse direction
and so the bundle $K_\gamma$ just becomes a scalar. Furthermore
$\zeta_{(\gamma;dz)}$ becomes the standard angle function for the
Lagrangian loops in $\C^n$, and hence Definition 2.6 reduces to
the standard Maslov index of the loop of Lagrangian Grassmanians.
\item
In fact, the above discussion does not depend on the particular choice
of the standard complex structure $j$ on $\R^{2n}$ but can be repeated
verbatim for the complex structure on $E_\gamma = \C^n$
induced from any almost complex structure $J$ compatible to the standard
symplectic structure $\omega_0$. Using the fact that the set of
compatible almost complex structure is contractible and the degree
is an homotopy invariant, the coisotropic Maslov index
is an invariant depending only on the homotopy class of
the complex line bundle $K_{(J; \gamma)}$. We note that this
homotopy class is an invariant depending only on
the coisotropic loop $\gamma$. In this sense, the
index $\mu(\gamma;\zeta)$ defined above is a symplectic invariant.
\item
A complete parallel discussion can be carried out in the case of
Calabi-Yau manifolds as we will discuss in section 7 later.
\endroster
\enddefinition

From the definition, it is clear that
to discuss invariance property of
the index $\mu(\gamma;\zeta)$, one needs more than just
a homotopy of the loop $\gamma$ but also
to dictate a proper condition for the section $\zeta$ as well.
This deviates from the homotopy invariance of the classical Maslov
index of Lagrangian loops [Ar].

However it has the property of symplectic invariance which we now
describe. Let $\gamma \in \Gamma_k(\R^{2n},\omega_0)$ and
$\zeta \in \Gamma(K_\gamma)$ be a Maslov section.
For any loop $A: S^1 \to Sp(2n)$, it naturally induces the push-forward
complex structure $A_*j$ on $\R^{2n}$ and a push-forward coisotropic
loop $A\cdot \gamma: S^1 \to \Gamma(\R^{2n},\omega_0)$.
It induces a pair $(A \cdot \gamma;A_*\zeta)$ of $A_*\zeta \in
K_{(A_*j,A\cdot \gamma)}$.

\proclaim{Theorem 2.9} Let $\gamma \in \Gamma_k(\R^{2n},\omega_0)$
and $\zeta \in \Gamma(K_{(j,\gamma)})$. Let $A: S^1 \to Sp(2n)$ be
any loop of symplectic matrices. Then we have
$$
\mu(A\cdot \gamma;A_*\zeta) = \mu(\gamma;\zeta). \tag 2.7
$$
\endproclaim

\medskip

\head{\bf \S 3. Graded coisotropic submanifolds,
and the coisotropic Maslov index}
\endhead

In this section, we will introduce the notion of grading
on coisotropic submanifolds and define an index of
maps carrying a transverse Maslov section with it.

First we denote by $\FF$ the null foliation of the
coisotropic submanifold $Y \subset (X,\omega)$
and consider the leafwise normal bundle $N\FF$ and its dual
$N^*\FF$ respectively.

We note that $J$ preserves both $N\FF \subset TX|_Y$ and $N^*\FF
\cong (T\FF\oplus NY)^\perp \subset T^*X|_Y$
and so induces decomposition of the complexifications
$$
\align
N\FF\otimes \C
& = N^{1,0}\FF \oplus N^{0,1}\FF \\
N^*\FF \otimes \C
& = (N^*)^{1,0}\FF \oplus (N^*)^{0,1}\FF.
\endalign
$$
By definition, the transverse canonical bundle is defined by
$$
K_Y \to Y = \det(\Lambda^{1,0}(N^*\FF)).
$$
In general $K_Y^{\otimes 2} \to Y$ is not trivial.
However it is trivial if $Y$ is either Lagrangian or
if the total space $X$ satisfies $2c_1(X,\omega_0)=0$.

This leads us to the following notion.

\definition{Definition 3.1}
\roster
\item We call a coisotropic submanifold
$Y \subset (X,\omega)$ {\it gradable} if its transverse
Maslov bundle $K_Y^{\otimes 2}$ is trivial.
\item We call a pair $(Y,[\Delta])$
a {\it graded} coisotropic submanifold where
$\Delta \in \Gamma(K_Y^{\otimes 2})$ and $[\Delta]$ its homotopy class.
Then we call the $\Delta$ {\it
transverse Maslov charge} and its homotopy class $[\Delta]$ a {\it
grading} of $Y$.
\endroster
\enddefinition

Here we note that our definition of graded 
coisotropic submanifolds is manufactured so that an index for
a map with boundary lying on the given coisotropic submanifold
can be defined which is always the case for Lagrangian submanifolds.
We would like to emphasize that it is not a generalization of the graded
Lagrangian submanifold used in [Ko] or [Se]. For example, according
to our definition, any Lagrangian submanifold is gradable and canonically
graded. 

Obviously for given compatible almost complex structure $J$
on $(X,\omega)$, the space of transverse Maslov charges of $Y$
is a principal homogeneous space of
$C^\infty(Y,S^1)$.  In other words, we have
$$
\Delta - \Delta' \in C^\infty(Y,S^1). \tag 3.1
$$

We are now ready to define the {\it coisotropic Maslov index} of a
map $w:(\Sigma,\del \Sigma) \to (X,Y)$ in the presence of
transverse Maslov charge $\Delta$ of $Y$ which we denote by $\mu_{(Y,
\Delta)}(w)$. We will suppress $\Delta$ from the notation,
whenever there is no danger of confusion.

Consider the pull-back $w^*TX$ and its symplectic
trivialization
$$
\Phi: w^*TX \to \Sigma \times \C^n.
$$
Such a trivialization, even the unitary one, always exists as long
as $\partial \Sigma \neq \emptyset$. Furthermore it is unique up
to homotopy when $\Sigma$ is a disc. We will restrict to the case
of discs in this paper and postpone more general discussions of
the higher genus cases elsewhere.  From now on, we assume $\Sigma$
is a disc $D = D^2$. We will denote by
$$
\{u_1, \cdots, u_n\}
$$
the corresponding unitary frame such that
$$
\text{span}_\C \{u_1, \cdots, u_k\}
= (TY^\omega\oplus JTY^\omega)^{\perp,J}.
\tag 3.2
$$

In the presence of transverse Maslov charge $\Delta$ on $Y$, we require that
the push-forward grading
$$
[\Delta_\Phi] : = [(\Phi_{\del D})_*(\Delta)]
$$
on the loop $\alpha_\Phi \in \Gamma_k(\R^{2n},\omega_0)$
coincides with the canonical grading
defined in section 2. In terms of the above unitary frame,
this means that
$$
[(\gamma)^*\Delta] = [(\gamma)^*(\theta^1 \wedge \cdots\wedge \theta^k)]
\tag 3.3
$$
where $\gamma= \del w: \del D \to Y$. More precisely, we have

\proclaim{Lemma 3.2} \roster
\item We can always choose a unitary frame
of $\gamma^*(TX)$ $\{u_1, \cdots, u_n\}$ over $\del D$
that satisfies (3.2) and (3.3).
\item We can then extend the above frame over $\del D$ to
that of $w^*TX$ over $D$.
\endroster
\endproclaim
\demo{Proof}
Recall that any complex vector bundle over $\del D \cong S^1$ is
trivial. This proves (1).
The statement (2) follows from the fact that
$\pi_2(U(n), U(k)\otimes\{id\}) = 0$ when $n\geq 2$ and $0\leq k\leq n$.
\qed\enddemo

We will call such a trivialization $\Phi$ an {\it admissible}
trivialization of $w^*(TX)$ with respect to $(Y,\Delta)$. We will
sometimes denote by $\Phi_\Delta$ if necessary to explicitly state
this choice.

Now we are ready to define the {\it coisotropic Maslov index
for the map $w:(D,\del D) \to (X,Y)$ with respect to the
transverse Maslov charge $\Delta$ on $Y$.

\definition{Definition 3.3}
Let $w$ be as above, and let $\Phi: w^*TX \to \Sigma \times \C^n$ be an
admissible trivialization. Then we define the index
$$
\mu_{(Y,\Delta)}(w) = [\Delta] - [\Phi^*(\zeta_{(\alpha_\Phi;dz)})] \in \Z
$$
and call it the coisotropic Maslov index of the map $w$
with respect to $\Delta$.
\enddefinition
It is easy to check the above definition does not depend on the
choice of admissible trivializations $\Phi_\Delta$.

In general, the coisotropic Maslov index will not be invariant under the
homotopy of maps $w_t: (D,\del D) \to (X,Y)$ for $t \in [0,1]$.
However this is so when the transverse Maslov charge $\Delta$ is also
{\it parallel}.
We will prove the following theorem in the next section.

\proclaim{Theorem 3.4} Suppose that the transverse Maslov charge
$\Delta$ on $Y$ is parallel with respect to the
connection induced from the canonical connection from $(X,\omega,J)$.
Let $\{(w_t)\}_{0 \leq t \leq 1}$ be a
smooth family of maps with boundary lying on $Y$. Then we have
$$
\mu_{(Y,\Delta)}(w_0) = \mu_{(Y,\Delta)}(w_1).
$$
\endproclaim

An obvious necessary condition for the existence of such a
parallel section is vanishing of the curvature of $K_Y$. This
theorem leads us to the following definition.

\definition{Definition 3.5}
We call a coisotropic submanifold $Y$ a {\it transversally
Ricci-flat} if $K_Y$ is flat.
\enddefinition

We will study various geometric properties of transversally
Ricci-flat coisotropic submanifolds in K\"ahler manifolds later in
section 5-7.

\head \bf \S 4. Coisotropic submanifolds in almost K\"ahler
manifolds
\endhead

In this section we recall the basic facts on the connection on the
almost K\"ahler manifolds $(X,\omega, J)$ following [Kb].

We now consider an almost complex structure $J$ that is compatible
to $\omega$, i.e, the triple $(X, \omega, J)$ defines an {\it
almost K\"ahler manifold} with the standard convention
$$
g = \omega(\cdot, J\cdot).
$$
We will choose a connection $\nabla$ on $TX$ which preserves both
$J$ and $\omega$. Such a connection always exist and in addition
becomes unique if we impose the condition that
$$
\nabla \omega =0 =\nabla J
$$
and the corresponding torsion form $\Theta$ is of type (0,2).

\definition{Definition 4.1}
A connection $\nabla$ of
the almost K\"ahler manifold $(X,\omega,J)$ is called the {\it
canonical} connection (see [Kb]), if its torsion form is
of type (0,2).
\enddefinition

\proclaim{Theorem 4.2 [Theorem 5.1, Kb]} Every almost K\"ahler
manifold $(X,\omega,J)$ carries the unique canonical connection.
This connection also satisfies
$$
\sum \Theta^i \wedge \bar\theta^i = 0. \tag 4.1
$$
\endproclaim

When we say that $(X,\omega,J)$ is an almost K\"ahler manifold,
we will always assume that it carries the canonical connection.
Now we need to do some basic calculations
involving the moving frame {\it adapted to} the given coisotropic
submanifold $Y \subset X$.

We choose an orthonormal frame of $X$
$$
\{e_1,e_2, \cdots, e_n, f_1, \cdots, f_n\} \tag 4.2
$$
such that
$$
\operatorname{span}_\R\{e_1, \cdots, e_k, e_{k+1}\cdots, e_n, f_1, \cdots,
f_k\}= TY \subset TX
$$
and
$$
\operatorname{span}_\R \{e_{k+1}, \cdots, e_n\} = (TY)^\omega : = E
\subset TY.
$$
The vectors
$$
u_j = \frac{1}{2}(e_j - if_j), \quad j = 1, \cdots, n
$$
form a unitary frame of $T^{(1,0)}X$ and the vectors
$$
\overline u_j = \frac{1}{2}(e_j + if_j), \quad j = 1, \cdots, n
$$
form a unitary frame of $T^{(0,1)}X$. Let
$$
\{e^*_1, \cdots, e^*_n, f^*_1, \cdots, f^*_n\}
$$
be the dual frame of (4.2).  The complex valued one forms
$$
\theta^j = \alpha^j + i \beta^j, \quad j = 1, \cdots, n
$$
form a unitary frame of $TX$ which are dual to the unitary frame
$\{u_1, \cdots, u_n\}$.

In terms of the above mentioned canonical connection $\nabla$, we
have
$$
\nabla u_j = \sum_i \omega^i_j \otimes u_i + \sum_\ell \tau_{\bar
j}^\ell \otimes \overline u_\ell \tag 4.3
$$
where $\omega_j^i$ is the connection one form
$$
\omega_j^i = \sum\omega_{jk}^i \theta^k + \omega_{j\bar k}^i
\bar\theta^k.
$$
The first structure equation with respect to the frame
$\{\theta^1, \cdots, \theta^n\}$ becomes
$$
d\theta^j = - \sum_i \omega^j_i \wedge \theta^i + \Theta^j \tag
4.4
$$
where the torsion form $\Theta^j$ is of the type (0,2)
$$
\Theta^j = \sum_i \tau^j_{\bar i} \wedge \overline \theta^i =
\sum_{i,\ell} N^j_{\bar{\ell} \bar{i}} \overline\theta^\ell\wedge
\overline\theta^i. \tag 4.5
$$
Furthermore if we set $N^j_{\bar{\ell}\bar i} = N_{\bar
j\bar{\ell}\bar i}$, then (4.1) implies
$$
N_{\bar j\bar{\ell}\bar i}= N_{\bar{\ell}\bar i\bar j}= N_{\bar
i\bar j\bar{\ell}}
$$
(see [Kb]), or equivalently
$$
\tau^j_{\bar i} = \tau^i_{\bar j}.
\tag 4.6
$$
The 2-nd structure equation of $\{\theta_1, \cdots, \theta_n\}$
is
$$
d\omega = -\omega\wedge \omega + K \tag 4.7
$$
where $K$ is the curvature 2-form.

The following proposition will play an important role in the proof
Theorem 3.4 for the maps $w:(D^2, \partial D^2) \to (X,Y)$
and also in our calculation of covariant derivative
of the transverse Maslov section later.

\proclaim{Proposition 4.3} Suppose that $J$ and $\nabla$ are as
above. Let $\gamma: [0,1] \to Y \subset X$ be a smooth curve on
$Y$ and let $\Pi_\nabla$ be the parallel translation from
$T_{\gamma(0)}X$ to $T_{\gamma(1)}X$ in $X$. Then $\Pi_\gamma$
maps $(E\oplus JE)|_{\gamma(0)} \subset T_{\gamma(0)}X$ to
$(E\oplus JE)|_{\gamma(1)} \subset T_{\gamma(1)}X$. In particular,
the subbundle $H_E:= E^{\perp,J} \subset TY$ is also invariant
under the parallel translation along such curves.
\endproclaim
\demo{Proof}
 Note that in terms of the metric $g$, it is enough to
prove that for any vector field $\eta$ on $Y$ such that $\langle
\eta(\gamma(0)), \xi_0 \rangle_g = 0$, the parallel translation
$\Pi_\gamma(\xi_0)$ also satisfies
$$
\langle \eta(\gamma(1)), \Pi_\gamma(\xi_0) \rangle = 0 \tag 4.8
$$
when $\xi_0 \in E\oplus JE$. We choose an orthonormal frame of $X$
$$
\{e_1, \cdots, e_n, f_1, \cdots, f_n\} \tag 4.9
$$
adapted to $Y$. The vectors
$$
u_j = \frac{1}{2}(e_j - if_j), \quad j = 1, \cdots, n
$$
form a unitary frame of $T^{(1,0)}X$ and the vectors
$$
\overline u_j = \frac{1}{2}(e_j + if_j), \quad j = 1, \cdots, n
$$
form a unitary frame of $T^{(0,1)}X$. Let
$$
\{\alpha_1, \cdots, \alpha_n, \beta_1, \cdots, \beta_n\}
$$
be the dual frame of (4.9)  The complex valued one forms
$$
\theta^j = \alpha^j + i \beta^j, \quad j = 1, \cdots, n
$$
form a unitary frame of $TX$ which is dual to the unitary frame
$\{u_1, \cdots, u_n\}$.

Substituting  (4.5) and
$$
\omega^j_i = \sum_\ell \omega^j_{i\ell} \theta^\ell +\sum_\ell
\omega^j_{i\bar\ell} \overline{\theta}^\ell
$$
into (4.4), we can write
$$
d\theta^j = -\sum_{i,\ell} \omega^j_{i\ell} \theta^\ell \wedge
\theta^i - \sum_{i,\ell} \omega^j_{i\bar\ell}
\overline{\theta}^\ell \wedge \theta^i + \sum_{i,\ell}
N^j_{\bar{\ell} \bar{i}} \overline\theta^\ell\wedge
\overline{\theta^i}.
$$
It is easy to see that for $d\theta^j$ with $j = 1, \cdots, k$ to
be in the ideal generated by $\{\theta^1, \cdots, \theta^k,
\overline\theta^1, \cdots, \theta^k\}$ on $Y$, we must have
$$
\omega^j_{i\ell} = \omega^j_{i \bar\ell} = N^j_{\bar i \bar \ell}=
0, \quad \text{for }\, i, \, \ell \geq k +1. \tag 4.10
$$
Let $\xi:[0,1] \to T^{(1,0)}X$ be the unique solution for
$$
\cases
\nabla_{\gamma'(t)} \xi = 0 \\
\xi(0) = \xi_0 \in E\oplus J E|_{\gamma(0)}
\endcases
$$
for given $\xi_0$. By definition of the parallel translation, we
have $\Pi_\gamma(\xi_0) = \xi(1)$. We compute
$$
\frac{d}{dt} \langle u_j(\gamma(t)), \xi(t) \rangle
 =  \langle \nabla_t u_j, \xi \rangle + \langle u_j, \nabla_t\xi
\rangle = \langle \nabla_t u_j, \xi \rangle.
$$
On the other hand we derive from (4.3) and (4.10),
$$
\nabla_t u_j = \sum_i^k \omega^i_j(\gamma'(t)) u_i + \sum_\ell^k
\tau_{\bar j}^\ell(\gamma'(t)) \overline u_\ell
$$
and hence
$$
\frac{d}{dt} \langle u_j(\gamma(t)), \xi(t) \rangle = \sum_i^k
\omega^i_j(\gamma'(t)) \langle u_i, \xi\rangle  + \sum_\ell^k
\tau_{\bar j}^\ell(\gamma'(t)) \langle \overline u_\ell, \xi
\rangle
$$
for $ 1 \leq i, \,  j \leq k$. Similarly we also derive the
equation for $\langle\overline u_i, \xi \rangle$ with $ 1 \leq i,
\, j \leq k$. Together we have a system of linear first order ODE
for $\langle u_i, \xi\rangle$ and $\langle \overline u_i, \xi
\rangle$ for $1 \leq i \leq k$ with the initial condition
$$
\langle u_i (\gamma(0), \xi(0) \rangle = \langle \overline
u_i(\gamma(0), \xi(0) \rangle  = 0.
$$
This proves that $\langle u_i (\gamma(t)), \xi(t) \rangle =
\langle \overline u_i(\gamma(t)), \xi(t) \rangle  \equiv 0$ for
all $t \in [0,1]$ and in particular at $t = 1$. This finishes the
proof. \qed\enddemo

\demo{Proof of Theorem 3.4} Because of the symplectic invariance
of the index $\mu(\gamma;\zeta)$, we may choose any trivialization
of $w_t^*TX$ for the definition of $\mu_{(Y,\Delta)}(w_t)$.
We denote the parameterized map
$$
W: [0,1] \times D \to X; \quad W(t,z) = w_t(z)
$$
and fix a trivialization of $W^*(TX)$
$$
\Phi: W^*(TX) \to [0,1] \times D \times \C^n
$$
Under this trivialization, $(\del w_t)^*TY$ gives rise to a loop
of coisotropic subspaces $\alpha_{w_t}: S^1 \to
\Gamma_k(T_{w(0)}X) \cong \Gamma_k(\C^n)$ for each $t \in [0,1]$.
Now we consider the parallel translations
$$
\Pi_t^0: (\del w_t)^*TX \to (\del w_0)^*TX
$$
along the paths
$$
t \mapsto w_t(\theta)
$$
for each $\theta \in \del D$.  The  pairs $(S_{\alpha_{w_t}},
L_{\alpha_{w_t}})$ corresponding to the loop $t \mapsto
\Pi_t^0\cdot \alpha_{w_t} \in \Gamma_k(\C^n)$ are mapped to a
$t$-parameter family of the pairs
$$
(\Pi^0_t(S_{\alpha_{w_t}}), \Pi^0_t(L_{\alpha_{w_t}}))
$$
in $\C^n$ and Proposition 4.3 implies that
$$
\Pi^0_t(S_{\alpha_{w_t}}) \equiv S_{\alpha_{w_0}} \subset \C^n.
\tag 4.11
$$
for all $t \in [0,1]$. Furthermore since we assume that $\Delta$
is parallel, $\Delta$ is invariant under the parallel translation
$\Pi^0_t$, i.e., we have
$$
\Pi_t^0(\Delta_t) = \Delta_0
\tag 4.12
$$
where $\Delta_t = (\del w_t)^*\Delta$. Then $\mu_{(Y,\Delta)}(w_0)
= \mu_{(Y,\Delta)}(w_1)$ follows from (4.11), (4.12), Proposition
4.3 and homotopy invariance of the usual Maslov index of the
Lagrangian Grassmanians. \qed\enddemo

\head \bf\S 5. Leafwise mean curvature vector
\endhead

From now on, we use indices $i, \, j,\dots$ from 1 to $n$, $a, \,
b, \, c, \cdots$ from 1 to $k$ and $\alpha, \, \beta, \, \cdots$
from $k+1, \cdots, n$ and use the summation convention for the
repeated indices. We restrict the 1-st structure equation (4.4) to
$Y$. Since the distribution $E$ on $Y$ is integrable and
$f^*_\alpha = 0$ on $Y$, both $d\theta^a$ and $d\bar\theta^a$
should be contained in the ideal generated by $\{\theta^1, \cdots,
\theta^k, \overline\theta^1, \cdots, \theta^k\}$. The 1-st
structure equation provides
$$
d\theta^a = - \omega^a_i \wedge \theta^i + \Theta^a = - \omega^a_b
\wedge \theta^b - \omega^a_\alpha \wedge \theta^\alpha + \Theta^a
$$
on $Y$. We will come back to the {\it non-integrable} case in the
appendix and restrict our attention to the integrable case
for the rest of the paper.

From now on, we will assume $(X,\omega, J)$ is K\"ahler, i.e.,
$\Theta \equiv 0$. In this case, the above first structure
equation becomes
$$
d\theta^a = - \omega^a_i \wedge \theta^i = - \omega^a_b \wedge
\theta^b - \omega^a_\alpha \wedge \theta^\alpha \tag 5.1
$$
Since this is contained in the ideal generated by $\{\theta^1,
\cdots, \theta^k, \overline\theta^1, \cdots, \theta^k\}$, and
$$
\omega_\alpha^a \wedge \theta^\alpha = \omega_{\alpha j}^a
\theta^j \wedge \theta^\alpha + \omega_{\alpha\bar j}^a \bar
\theta^j \wedge \theta^\alpha
$$
we derive
$$
\omega_{\alpha\beta}^a \theta^\beta \wedge \theta^\alpha +
\omega_{\alpha\bar\beta}^a \bar\theta^\beta \wedge \theta^\alpha =
0
$$
and so
$$
\omega_{\alpha\beta}^a = \omega_{\beta\alpha}^a, \quad
\omega_{\alpha \bar\beta}^a = 0. \tag 5.2
$$
Now we need to incorporate some {\it real} aspect of the
submanifold $Y \subset (X,\omega)$. For this, we write the
connection one form $\omega^i_j$ as
$$
\omega^i_j = \alpha_j^i + i \beta_j^i
$$
where $\alpha_j^i, \, \beta_j^i$ are real one forms.
The unitarity of connection, i.e, the skew-Hermitian property of
$(\omega_j^i)$ implies
$$
\alpha_j^i = - \alpha_i^j, \quad \beta_j^i = \beta_i^j. \tag 5.3
$$
The first structure equation $d\theta^\ell = -\omega_j^\ell \wedge \theta^j$
becomes
$$
\cases
de^*_\ell = - \alpha_j^\ell \wedge e_j^* + \beta_j^\ell \wedge f_j^* \\
df^*_k = - \beta_j^k \wedge e_j^* - \alpha_j^k \wedge f_j^*
\endcases
\tag 5.4
$$
For $k+1 \leq  \alpha \leq n$, since $f^*_\alpha = 0$ on $Y$, we derive
$$
0 = df^*_\alpha = - \beta_j^\alpha \wedge e_j^* - \alpha_b^\alpha
\wedge f_b^* \tag 5.5
$$
on $Y$.  By Cartan's lemma, we conclude
$$
\aligned
\beta_j^\alpha & = A_{j\ell}^\alpha e_\ell^* + B_{jb}^\alpha f_b^* \\
\alpha_b^\alpha & = C_{bj}^\alpha e_j^* + D_{bc}^\alpha f_c^*
\endaligned
\tag 5.6
$$
where the coefficients satisfy
$$
A^\alpha_{j\ell} = A^\alpha_{\ell j}, \, B_{jb}^\alpha =
C_{bj}^\alpha, \, D_{bc}^\alpha = D_{cb}^\alpha. \tag 5.7
$$
Therefore the second fundamental form $S = \sum_\alpha S^\alpha f_\alpha$
of $Y$ is given by the symmetric matrix
$$
S^\alpha = \Big( \matrix A^\alpha_{\ell j}, & C^\alpha_{jb}  \\
B^\alpha_{bj}, & D^\alpha_{cb}
\endmatrix\Big).
$$
Now we are ready to define the {\it leafwise mean curvature vector}.

\definition{Definition 5.1} Let $Y \subset X$ be a coisotropic submanifold
and let $S = S^\alpha f_\alpha$ be the second fundamental form of
$Y \subset (X,g)$. The leafwise mean curvature vector of $Y$ is the
partial trace of $S$
$$
\vec H^\parallel: = S(e_\beta, e_\beta) =
S^\alpha(e_\beta,e_\beta) f_\alpha. \tag 5.8
$$
The {\it leafwise mean curvature one form} of $Y$ is defined by
$$
\alpha_Y^\parallel : = \vec H^\parallel \rfloor \omega. \tag 5.9
$$
\enddefinition

It is easy to check that the right hand side of (5.8) is
independent of the choice of the frame adapted to $Y$ and so the
leafwise mean curvature vector and so the leafwise mean curvature one form
are well-defined.

In the above moving frame, the leafwise mean curvature vector is given
by
$$
\vec H^\parallel = A^\alpha_{\beta\beta} f_\alpha
$$
Recalling the K\"ahler form is given by
$\omega = \frac{i}{2} \theta^j \wedge \bar \theta^j = e_k^* \wedge f_k^*$,
we prove that the mean curvature one form becomes
$$
\vec H^\parallel \rfloor \omega =
- A^\alpha_{\beta\beta} e^*_\alpha \tag 5.10
$$

We summarize the above calculations into the following proposition

\proclaim{Proposition 5.2} Let $Y \subset (X,\omega)$ be
coisotropic and $(X,\omega,J)$ be K\"ahler. Then the leafwise mean
curvature one form has the formula
$$
\alpha_Y^\parallel:= \vec H^\parallel \rfloor \omega|_{TY} =
- A^\beta_{\alpha\alpha}e_\beta^* = - A^\alpha_{\beta\alpha}e_\beta^*.
\tag 5.11
$$
\endproclaim

Now we will derive leafwise derivatives of various leafwise
differential forms, including $\alpha_Y^\parallel$. To carry out
this derivation in a coherent manner, we need to provide a brief
review  the concept of $E$-de Rham complex [NT] associated to the
structure of Lie algebroid [Mac].

\definition{Definition 5.3} Let $M$ be a smooth manifold.
A {\it Lie algebroid} on $M$ is a triple $(E,\rho, [\, ,\, ])$, where
$E$ is a vector bundle on $M$, $[\, ,\, ]$ is a Lie algebra structure
on the sheaf of sections of $E$, and $\rho$ is a bundle map
$$
\rho: E \to TM \tag 5.12
$$
such that the induced map
$$
\Gamma(\rho): \Gamma(M;E) \to \Gamma(TM) \tag 5.13
$$
is a Lie algebra homomorphism and, for any sections $\sigma$ and $\tau$
of $E$ and a smooth function $f$ on $M$, the identity
$$
[\sigma, f\tau] = \rho(\sigma)[f]\cdot \tau + f\cdot [\sigma,
\tau]. \tag 5.14
$$
\enddefinition

\definition{Definition 5.4 [Definition 2.2, NT]}
Let $(E,\rho,[\, ,\, ])$ be a Lie algebroid on $M$. The $E$-de Rham complex
$(^E\Omega^\bullet(M), ^Ed)$ is defined by
$$
\aligned
& ^E\Omega(\Lambda^\bullet(E^*)) = \Gamma(\Lambda^\bullet(E^*))\\
& ^Ed\omega(\sigma_1, \cdots, \sigma_{k+1}) = \\
&\quad = \sum_i(-1)^i \rho(\sigma_i)\omega(\sigma_1, \cdots, \widehat{\sigma_i},
\cdots, \sigma_{k+1}) \\
&\qquad + \sum_{i < j}(-1)^{i+j-1}\omega([\sigma_i,\sigma_j],\sigma_1,
\cdots,\check\sigma_i, \cdots, \check\sigma_j, \cdots, \sigma_{k+1}).
\endaligned
\tag 5.15
$$
The cohomology of this complex will be denoted by $^EH^*(M)$ and called
the $E$-de Rham cohomology of $M$.
An {\it $E$-connection} on a vector bundle $F$ on $M$ is a linear map
$$
\nabla: \Gamma(F\otimes \Lambda^\bullet(E^*)) \to \Gamma(F\otimes
\Lambda^{\bullet + 1}(E^*)) \tag 5.16
$$
satisfying the Leibnitz rule:
$$
\nabla(f\sigma) = ^Ed(f)\cdot \sigma + f\cdot \nabla\sigma.
$$
Similarly we can define the notion of {\it $E$-curvature}
$^E K_\nabla \in \Gamma(\Lambda^2(E^*) \otimes \text{End}(F))$
in an obvious way.
\enddefinition
We would like to emphasize that $^EH^*(M)$ is not a topological
invariant of $Y$ but an invariant of the Lie algebroid $E =
TY^\omega$ or of the null foliation $\FF$ of $Y$. For example, the
class is not invariant under the general homotopy but invariant
only under the homotopy {\it tangential} to the foliation in an
obvious sense. We can also denote $TY^\omega = T\FF$ the tangent
bundle of the null foliation $\FF$ and by $E^* = T^*\FF$ the
cotangent bundle of $\FF$.

In our case, $M = Y$ and $E: = TY^\omega$ and the anchor map
$\rho: E \to TY$ is nothing but the inclusion map $i: TY^\omega
\to TY$. The integrability of $TY^\omega$ implies that the
restriction of the Lie bracket on $\Gamma(TY)$ to
$\Gamma(TY^\omega)$ defines the Lie bracket $[\, ,\, ]$ on
$\Gamma(E)$. Therefore the triple
$$
(E=TY^\omega, \rho = i, [\, ,\, ])
$$
defines the structure of Lie algebroid and hence the
$E$-differential and $E$-connections. In our case, the
corresponding $E$-differential is nothing but $d_\FF$ the exterior
derivative along the null foliation $\FF$. and its cohomology,
denoted by  $H^*(Y,\omega)$, the cohomology $H^*(\FF)$ of the
foliation $\FF$. We denote the $E$-connection of a vector bundle
$F$ over $Y$ of this Lie algebroid by $\nabla^\omega$ in general.

Now we go back to further study the leafwise mean curvature vector.
The the Ricci form $i \rho$ of the metric $g
= \omega(\cdot, J\cdot)$ represents the curvature of the canonical
line bundle
$$
K : = \Lambda^{n,0}(T^*X \otimes \C).
$$
We first recall one geometric fact about Lagrangian submanifolds
on the general K\"ahler manifolds.  For a given Lagrangian
embedding $i:L \hookrightarrow (X, \omega)$, we consider the
pull-back bundle $F=i^*K \to L$. If $L$ is orientable, the line
bundle $i^*K$ is always trivial. In general, the square
$i^*K^{2\otimes}$ is trivial. In fact, for the orientable
Lagrangian submanifold, there exists a {\it canonical} section
which we denote by $\Omega_L$: we choose a positively oriented
Darboux frame
$$
\{e_1, \cdots, e_n, f_1, \cdots, f_n\}
$$
so that $\{e_1, \cdots, e_n\}$ spans, $TL$ and the associated
holomorphic frame $\{\theta_1, \cdots, \theta_n\}$ of
$\Lambda^{n,0}(TX\otimes \C)$ adapted to $L$. In particular, the
$n$-form
$$
e_1^* \wedge \cdots \wedge e_{n}^*
$$
provides the volume form of $Y$ with respect to the induced
metric. Then the complex $n$-form
$$
\Omega_L = \theta_1 \wedge \cdots \wedge \theta_n
$$
where $\theta_j = e_i^* + i f_i^*$ provides the canonical section.

Now we study the general coisotropic cases. The discussion below
will be parallel with the case in $(\R^{2n},\omega_0)$ considered
in section 2.  Recalling that the direct sum
$$
\EE:= T\FF \oplus N_JY \subset TX|_Y
$$
is invariant under the
action of $J : TX|_Y \to TX|_Y$, we consider the complexification
$$
\EE \otimes \C =  T^{1,0}\EE \oplus T^{0,1}\EE
$$
where the splitting is with respect to the
complex structure $J: T\FF \oplus N_JY \to T\FF \oplus N_JY$.
One can also consider the dual bundle
$$
\EE^*= N^*\FF \oplus T^*\FF \cong N_JY\oplus T\FF =
\EE. \tag 5.17
$$
Then
we consider the top exterior power of $\EE^{1,0}$,
$\det \EE^{1,0}$, as a complex line vector bundle.

Note that the complex line bundle
$\det\EE^{1,0})^{\otimes 2}$ has a canonical section given by
$$
(u_{k+1} \wedge \cdots \wedge u_n)^{\otimes 2}
$$
for any adapted frame $\{e_1, \cdots, e_n ,f_1, \cdots, f_n\}$ and its
associated Hermitian frame.  It is easy to check that this
form does not depend on the choice of the adapted frame (4.9) and
so globally well-defined.

We will now compute the covariant derivative
$$
\nabla (u_{k+1} \wedge \cdots \wedge u_n)
$$
with respect to the induced connection on $(\det\EE^{1,0})^{\otimes 2}$
from that of $(T\FF \oplus J(T\FF)) \otimes \C$. More precisely
we have the following

\proclaim{Lemma 5.5} The canonical connection on $TX|_Y$
naturally induces a connection on $T\FF \otimes J (T\FF) \subset TX|_Y$
by restriction.
\endproclaim
\demo{Proof} This is an immediate consequence of Proposition 4.3.
\enddemo

We will also compute the $E$-covariant derivative
$$
\nabla^\omega (u_{k+1} \wedge \cdots \wedge u_n)
$$
where the $E$-connection $\nabla^\omega$ is nothing but the
restriction of the induced connection on $F = \det \EE^{1,0}$ from
the connection $\nabla$ on $T\FF \otimes J(T\FF)$ defined by Lemma
5.5. One can easily check that $\nabla^\omega$ really satisfies
the defining property $E$-connection in Definition 5.4.

\proclaim{Proposition 5.6} Let $(X,\omega, J)$ be K\"ahler and $i:
Y \hookrightarrow (X,\omega)$ be a coisotropic embedding. We
denote by $\nabla$ the connection on $\det \EE^{1,0} \to Y$
defined in Lemma 5.5.  Then we have
$$
\nabla (u_{k+1} \wedge \cdots \wedge u_n) =
i^*(\omega^\alpha_\alpha)\cdot (u_{k+1} \wedge \cdots \wedge u_n)
\tag 5.18
$$
and
$$
\nabla^\omega (u_{k+1} \wedge \cdots \wedge u_n) = -i
\alpha_Y^\parallel \cdot (u_{k+1} \wedge \cdots \wedge u_n) \tag
5.19
$$
In particular, we have
$$
curv(\det(\EE^{1,0})) = i^*(d\omega_\alpha^\alpha) \tag 5.20
$$
for the curvature of the line bundle $\det(\EE^{1,0})$ is given by
and its $E$-curvature is given by
$$
^Ecurv(\det(\EE^{1,0})) = - i d_\FF(\alpha_Y^\parallel). \tag 5.21
$$
\endproclaim

As in the Lagrangian case [Oh2], the following is an immediate
corollary of Proposition 5.6  whose proof we omitted.

\proclaim{Corollary 5.7 (Compare with [Corollary 3.3, Oh2])} Let
$i:Y \to X$ be any coisotropic embedding. Then the
holonomy of the bundle $\det(\EE^{1,0})$ over a loop $\gamma \subset Y$ with
respect to the induced connection is given by
$$
\exp(-i \int_\gamma \alpha_Y^\parallel)
$$
provided that $\gamma$ is tangent to the null foliation of $Y$,
i.e., $\gamma'(t) \in TY^\omega$ for all $t$.
\endproclaim

We first state the following symmetry property of the second
fundamental form $B$ of coisotropic submanifolds $Y$ in a K\"ahler
manifolds, which is the analog of the symmetry property for the
Lagrangian submanifolds [Lemma 3.10,Oh1].

\proclaim{Lemma 5.8} Let $B$ be the second fundamental form of the
coisotropic submanifold $Y \subset X$. Consider the tri-linear
form on $T\FF$ defined by
$$
(X,Y,Z) \mapsto \langle B(X,Y), JZ \rangle.
$$
Then we have
$$
\langle B(X,Y), JZ \rangle = \langle B(X,Z), JY \rangle
$$
for all $X, \, Y, \, Z \in T\FF$.
\endproclaim
\demo{Proof} This is a re-statement of the property that
$A^\alpha_{\beta\gamma}$ is fully symmetric over $\alpha,\, \beta$
and $\gamma$ which follows from (5.3) and (5.7). \qed\enddemo

\demo{Proof of Proposition 5.6} For any $X \in TY$, we compute
$\nabla_X(u^{k+1} \wedge \cdots \wedge u^{n})$ with $u_{k+1}
\wedge \cdots \wedge u_{n}$ taken as a unit frame of $\det
\EE^{1,0}$. Then we have
$$
\nabla_X(u_{k+1} \wedge \cdots \wedge u_{n}) = \sum_{\alpha=k+1}^n
u_{k+1} \wedge \cdots \wedge \nabla_X u_\alpha \wedge \cdots
\wedge u^n. \tag 5.22
$$
From the way how the connection $\nabla$ on $T\FF \oplus J(T\FF)$ is
defined in Lemma 5.5, we have
On the other hand Proposition 4.3 together with the structure equations
(4.1) restricted to $Y$ implies that the covariant derivative becomes
$$
\nabla_X u_\alpha  = \omega_\alpha^\gamma(X)u_\gamma
$$
and hence
$$
\nabla_X (u_{k+1} \wedge \cdots \wedge u_n) =
\omega^\alpha_\alpha(X)(u_{k+1} \wedge \cdots \wedge u_n).
\tag 5.23
$$
This proves (5.20). On the other hand by restricting to $T\FF$ and
using the fact that $\alpha^\gamma_\gamma = 0$ we have
$$
\nabla_X (u_{k+1} \wedge \cdots \wedge u_n) =  i
\beta^\alpha_\alpha(X)(u_{k+1} \wedge \cdots \wedge u_n).
\tag 5.24
$$
for $X \in T\FF$.
On the other hand, we have already shown that $\beta^\gamma_\gamma
= A^\gamma_{\gamma\alpha} e_\alpha^* = -\alpha_Y^\parallel$ in
(5.6), (5.11) and the symmetry property of
$A^\alpha_{\beta\gamma}$. This finishes the proof (5.21).
\qed\enddemo

\head{\bf \S 6. Pre-K\"ahler manifolds and its transverse
canonical bundle}
\endhead

The beginning discussion in this section is intrinsic in that it depends
only on the corresponding {\it pre-symplectic structure}
$\omega_Y$ on $Y$, while the discussions in the previous sections
are extrinsic in that it describes the property of the coisotropic
embedding into $(X,\omega)$. In the end of the section, we will
derive a compatibility condition between them.

For this purpose, it seems to useful to define another intrinsic notion

\definition{Definition 6.1} Let $(Y,\omega)$ be a pre-symplectic manifold
and fix a projection $\Pi:TY \to TY$ and its associated splitting
$TY = T\FF \oplus G_\Pi$ Let $J$ be an almost complex structure on
the normal bundle $N\FF \subset TY$ compatible to the
pre-symplectic structure $\omega_Y$ in the sense that
$$
g|_{N\FF} = \omega_Y(\cdot, J\cdot)|_{N\FF}
$$
We call the triple $(Y,\omega_Y, J)$ a {\it pre-K\"ahler}
manifold. We denote by $curv(K_Y)$ and
$^Ecurv(K_Y)$ the curvature and $E$-curvature of the
transverse canonical bundle $K_Y$ of $Y$, respectively.
\enddefinition

Obviously any coisotropic submanifold $Y$
in a K\"ahler manifold $(X,\omega, J)$
carries the induced pre-K\"ahler structure.

Considering $\theta^1\wedge \cdots \wedge \theta^k$ as a local
frame of $K_Y$, a straightforward computation shows that
the covariant derivative
$$
\nabla_X(\theta^1 \wedge \cdots \wedge \theta^k) =-
\omega^a_a(X)(\theta^1 \wedge \cdots \wedge \theta^k) \tag 6.1
$$
for any $X \in TY$. Here again we use Proposition 4.3 to
define the natural connection on the bundle $H_E:= (T\FF)^{\perp,J}$
induced from $TX|_Y$. The calculation is done with this
natural connection.

Therefore the curvature of $K_Y$ is given by
$$
curv(K_Y) = - i^*(d\omega^a_a). \tag 6.2
$$
On the other hand if we denote by
$$
i^*_\FF(\omega_a^a)
$$
the pull-back of the form $\omega^a_a$ to the leaves of the
null-foliation $\FF$,
The the $E$-curvature of the bundle $K_Y$ is given by the leafwise
differential
$$
^Ecurv(K_Y) = - d_\FF(i^*_\FF\omega_a^a). \tag 6.3
$$
These forms do not depend on the choice of frames adapted to $Y$
but depends only on the triple $(Y,\omega_Y,J_Y)$.

We now derive the following which is the coisotropic analog of the
result by Morvan [Mo] and Dazord [D], which relates the intrinsic
curvatures $(Y,\omega_Y,J_Y$
and the extrinsic curvatures of $Y \subset (X,\omega, J)$
and the ambient Ricci-curvature of $(X,\omega,J)$.

\proclaim{Theorem 6.2} Let $(X,\omega, J)$ K\"ahler and $i: Y
\subset (X,\omega,J)$ be any coisotropic submanifold. Let $K =
(K^i_j)$ be the curvature two form and $i \rho = K^j_j$ the Ricci
form of the metric. Then we have the formula
$$
\align
-curv(K_Y) & = - curv(\det\EE^{1,0}) + i i^*(\rho) \tag 6.4 \\
- ^Ecurv(K_Y) & = i d_\FF(\alpha_Y^\parallel) + i i^*_\FF(\rho)
\tag 6.5
\endalign
$$
where the two form $i^*_\FF(\rho)$ is the restriction to
$TY^\omega=T\FF \subset TY$ of the Ricci form $\rho$ of $X$.
Equivalently, we have
$$
\align
& - curv(K_Y) + curv(\det\EE^{1,0}) =i i^*(\rho) \tag 6.6 \\
& - ^Ecurv(K_Y)  - i d_\FF(\alpha_Y^\parallel) =i i^*_\FF(\rho)
\tag 6.7
\endalign
$$
\endproclaim
\demo{Proof} We will follow the above used notations in section 4.
We first note that
$$
d\beta^\beta_\gamma = \operatorname{Im}(d\omega^\beta_\gamma).
\tag 6.8
$$
Now we note that {\it if we restrict $\beta^\alpha_\gamma$ in
(5.6) to the leaves of $E$}, we have
$$
i^*_\FF(\beta^\alpha_\gamma) = A^\alpha_{\gamma\mu}e_\mu^*.
$$
On the other hand, taking the trace of the 2-nd structure equation
we have
$$
d\omega_j^j = -\omega^j_i\wedge \omega^i_j + K_j^j = K_j^j =
i\rho. \tag 6.9
$$
We decompose the left hand side of (6.9) and rewrite
$$
d\omega^\alpha_\alpha + d\omega_a^a = i\rho.
$$
This immediately proves  (6.6).
Restricting this to the leaves of $\FF$, we have
$$
d_\FF(i^*_\FF(\omega^\alpha_\alpha)) + d_\FF(i^*_\FF\omega_a^a) =
i i^*_\FF\rho. \tag 6.10
$$
On the other hand, we derive from (6.8)
$$
d_\FF(i^*_\FF(\omega_\beta^\beta)) = i
d_\FF(i^*_\FF(\beta_\beta^\beta)) = -i d_\FF(\alpha_Y^\parallel)
\tag 6.11
$$
and
$$
-^Ecurv(K_Y) =  d_\FF(i^*_\FF(\omega_a^a))
$$
from the definition (6.3) of the $E$-curvature of $K_Y$. This
finishes the proof. \qed\enddemo

\proclaim{Corollary 6.3} If $(X,\omega,J)$ is Einstein-K\"ahler,
i.e, $\rho = \lambda \omega$, then we have
$$
\align
& curv(K_Y) - curv(\det \EE^{1,0}) = 0. \tag 6.12 \\
& ^Ecurv(K_Y) + i d_\FF(\alpha_Y^\parallel) = 0. \tag 6.13
\endalign
$$
\endproclaim

This is the coisotropic analog to the well-known fact that the
mean curvature one form is always closed for the Lagrangian
submanifolds in Einstein-K\"ahler manifolds.

Now we restrict to the
case when $Y$ carries a parallel transverse Maslov charge $\Delta$.
An obvious necessary condition for $Y$ to carry such a
parallel transverse Maslov charge is the vanishing of its ordinary
curvature
$$
curv(K_Y) = 0 \tag 6.14
$$
i.e., $Y$ must be transversally Ricci-flat. Therefore we have
the following

\proclaim{Proposition 6.4} Suppose that $Y$
satisfies (6.14), i.e., is transversally Ricci-flat.  Then we have
$$
-id_\FF(\alpha_Y^\parallel) = i^*_\FF(\rho).
$$
In particular, for the Einstein-K\"ahler case, the mean leafwise
mean curvature form is leafwise closed and so define an
infinitesimal deformation of $Y$ as a coisotropic submanifold.
\endproclaim

\head \bf \S 7. Special coisotropic submanifolds in Calabi-Yau manifolds
\endhead

In this section, we restrict to the case when $(X,\omega, J)$ is
Calabi-Yau and study special features of geometry of coisotropic
submanifolds thereof. A coisotropic submanifold will have an induced
pre-K\"ahler structure $(Y,\omega_Y,J_Y)$ from $X$.

Let $\Omega$ be a holomorphic volume form of $X$ with constant
length one or equivalently a non zero holomorphic section of the
canonical line bundle $K$ with length one.

For any given coisotropic submanifold $Y \subset (X,\omega)$,
taking the restriction
of $\Omega$ to the leaves of $\FF$, we define a bundle map
$$
\widetilde\Omega_Y: (\det(T^{1,0}\EE))^{\otimes 2} \to
(K_Y)^{\otimes 2}
$$
by
$$
(\xi)^{\otimes 2} \mapsto (\xi \rfloor \Omega)^{\otimes 2}.
$$
In particular, this pushes forward the global section $(u_{k+1}
\wedge \cdots \wedge u_n)^{\otimes 2} \in \Gamma(\EE^{1,0})$ to
a global section on $K_Y^{\otimes 2}$. Therefore we have

\definition{Definition \& Proposition 7.1} Any coisotropic submanifold $Y$
in Calabi-Yau manifolds is canonically graded by the section
induced by $\Omega$. We denote this canonical transverse Maslov charge
by $\Delta_{(\Omega;Y)}$ and call it
the {\it canonical transverse Maslov charge of $Y$} and the
corresponding grading the {\it canonical grading}.
\enddefinition

This canonical grading induced by $\Omega$
enables us to define the coisotropic
Maslov index for the maps $w: (D,\del D) \to (X,Y)$
for all coisotropic submanifolds $Y$ in $X$ simultaneously, which
we denote by $\mu_{(Y;\Omega)}$.

Furthermore for the given frame adapted to $Y$, we can also write
$$
\Omega = f \theta^1 \wedge \cdots \wedge \theta^n
\tag 7.1
$$
for a locally defined function $f$ with values in $S^1$.
Since $\Omega$ is parallel, we have
$$
\align 0 & = \nabla_X(\Omega)= \nabla_X(f \theta^1 \wedge\cdots
\wedge \theta^n) \\
& = (df(X) - f \omega^j_j(X))\theta^1 \wedge\cdots \wedge \theta^n
\endalign
$$
and hence
$$
df(X) - f \omega^j_j(X) = 0
$$
for any $X \in TY$.  In the same frame, we have
$$
(u_{k+1} \wedge \cdots \wedge u_n) \rfloor \Omega
 = f (\theta^{1} \wedge \cdots \wedge
\theta^k). \tag 7.2
$$

Now we compute the covariant derivatives of both sides of (7.2)
separately.
From the right hand side, we derive
$$
\nabla_X (f \theta^1 \wedge \cdots\wedge \theta^k) = (df(X) - f
\omega^b_b(X)) (\theta^1 \wedge \cdots \wedge \theta^k) \tag 7.3
$$
for $X \in TY$. For the left hand side, using (5.18) and the fact
that $\Omega$ is parallel, we derive
$$
\nabla_X(u_{k+1} \wedge \cdots \wedge u_n) \rfloor \Omega)
= \omega^\alpha_\alpha(X)(u_{k+1} \wedge \cdots \wedge u_n) \rfloor \Omega).
\tag 7.4
$$
Comparing (7.3) and (7.4) with the identity (7.2), we derive
$$
f \omega^\alpha_\alpha(X) = df(X) - f \omega^a_a(X)
$$
for all $X \in TY$, i.e., we have
$$
0 = d \gamma - i^*(\omega^a_a) - i^*(\omega^\alpha_\alpha)
= d\gamma - i^*(\omega_j^j)
\tag 7.5
$$
where $\gamma = \log f$.
Restricting to $T\FF \subset TY$, we also have
$$
0 = d_\FF \gamma - i^*_\FF(\omega^a_a) + i\alpha_Y^\parallel
\tag 7.6
$$
Note that the obvious integrability conditions for the equation (7.5)
and (7.6) are
$$
d(i^*(\omega_a^a) +i^*(\omega_\alpha^\alpha)) = 0.
\tag 7.7
$$
and
$$
d_\FF(i^*_\FF(\omega_a^a) -i \alpha_Y^\parallel) = 0
\tag 7.8
$$
respectively, which, we have already shown in section 6,
holds on any Einstein-K\"ahler manifold
and so does on Calabi-Yau manifold.

\definition{Definition 7.2} Let $(X,\omega,J)$ be Calabi-Yau. We call
$Y \subset (X,\omega,J)$ {\it special coisotropic submanifold}
(respectively {\it leafwise special coisotropic submanifold}) if the
canonical section $\Delta_{(Y,\Omega)}$ is parallel (respectively
leafwise parallel), or equivalently if
$$
d\gamma - i^*(\omega_a^a) = 0. \tag 7.9
$$
and
$$
d_\FF\gamma - i^*_\FF(\omega_a^a) = 0 \tag 7.10
$$
with respect to the given frame, respectively.
\enddefinition
Obviously one necessary condition for the embedding $Y \subset (X,\omega,J)$
is special coisotropic is that the pre-K\"ahler structure $(Y,\omega_Y,J_Y)$
must be transversely Ricci-flat.

We recall that for the Lagrangian case, $N\FF = \{0\}$ and
$$
(i^*\Omega)^{\otimes 2} = g (\Omega_L)^{\otimes 2}
$$
for $g:L \to S^1$ which is a globally well-defined function and
called the {\it angle function} of the Lagrangian submanifold.
Therefore in this case, the above condition means
$$
g = (u_1\wedge \cdots \wedge u_n \rfloor \Omega)^{\otimes 2}
$$
satisfies $dg = 0$ on $L$, i.e., $f$ is constant.
Therefore Definition 7.2 reduces to
the usual special Lagrangian condition for the Lagrangian case.

Furthermore for the Lagrangian case, special Lagrangian
condition implies {\it minimality} of Lagrangian
submanifolds and vice versa (at least for orientable Lagrangian
submanifolds).

The following theorem is the coisotropic analog to this fact.

\proclaim{Theorem 7.3} Let $(X,\omega, J)$ be Calabi-Yau and $Y
\subset (X,\omega,J)$ be a leafwise special coisotropic submanifold.
Then $Y$ is leafwise special coisotropic if and only if it satisfies
$Y$ is leafwise minimal, i.e.,
$$
\alpha_Y^\parallel = 0
$$
\endproclaim
\demo{Proof}
This immediately follows from (7.5).
\qed\enddemo

Now we relate the above geometric study of
coisotropic submanifolds in Calabi-Yau with the coisotropic
Maslov index of the maps $w:(D,\del D) \to (X,Y)$.

For given $w$, we choose our unitary frame $\{v_1, \cdots, v_n\}$
on $w^*(TX)$ so that it is admissible to
the transverse Maslov charge $\Delta_{(Y;\Omega)}$.
We denote by $\{\phi^1, \cdots, \phi^n\}$ its dual frame.
The admissibility means
$$
(\del w)^*(\Delta_{(Y,\Omega)}) = (\phi^1 \wedge\cdots\wedge
\phi^k)^{\otimes 2}.
\tag 7.11
$$
Now we write
$$
(v_{k+1} \wedge \cdots \wedge v_n)^{\otimes 2}
= g\cdot (e_{k+1} \wedge \cdots \wedge e_n)^{\otimes 2}
$$
for $g: \del D \to S^1$, where we recall $T\FF = \text{span}_\R\{e_{k+1},
\cdots, e_n\}$.

The Maslov index $\mu_{(Y;\Omega)}(w)$ is computed by
the degree of the map
$$
g : \del D \to S^1
$$
$$
\mu_{(Y,\Omega)}(w) = \frac{1}{2\pi}\int_{\del D} g^*d\theta
\tag 7.12
$$
where $d\theta$ is the canonical angular form on $S^1$.

On the other hand, we can choose our frame
$\{v_1', \cdots, v_n'\}$ of $w^*TX$, which is
{\it not necessarily admissible}, so that
$$
v'_j = u_j \quad \text{on } \, (\del w)^*TY
\tag 7.13
$$
where $\{u_1, \cdots, u_n\}$ the restriction of
a frame of $(\del w)^*TX$ adapted to $(\del w)^*TY$ over $\del D$.
We write
$$
\Omega = f\cdot (\phi')^1 \wedge \cdots \wedge (\phi')^n
$$
for some $f : \del D \to S^1$. We write
$$
\Big(v'_{k+1} \wedge \cdots \wedge v'_n) \rfloor \Omega)^{\otimes 2} =
f^2 ((\phi')^1 \wedge \cdots \wedge (\phi')^k))^{\otimes 2}
\tag 7.14
$$
In particular, recalling Definition 3.3 of $\mu_{(Y;\Omega)}(w)$,
(7.14) implies
$$
\mu_{(Y;\Omega)}(w) = \text{deg}(f^2).
$$
Noting that the calculation in the beginning of this section is
purely on $Y$, we will have exactly the same equation as in
(7.5) with the connection one form $\omega^i_j$ replaced by
$(\omega')_j^i$ with respect to the frame $\{v_1',\cdots, v_n'\}$.
However when restricted to $Y$, we have
$$
(\omega')^i_j = \omega^i_j.
$$

This proves the following local index formula

\proclaim{Theorem 7.4} For any coisotropic submanifold $Y$ and
a frame $\{u_1, \cdots, u_n\}$ of $(\del w)^*TX$ adapted to
$(\del w)^*TY$, we have
$$
\mu_{(Y;\Omega)}(w)
= \frac{i}{\pi}\int_{\del w} i^*(\omega^j_j)
= \frac{i}{\pi}\int_{\del w} (i^*(\omega^a_a) + i^*(\omega^\alpha_\alpha))
\tag 7.15
$$
for all $w:(D,\del D) \to (X,Y)$. In particular, coisotropic Maslov
index depends only on the boundary map $\del w$.
\endproclaim

\proclaim{Corollary 7.5} Let $Y \subset (X,\omega)$ be a coisotropic
embedding. If the unitary
frame $\{u_1, \cdots, u_n\}$ of $TX|_Y$ adapted to $Y$ extends over the
neighborhood of the image of $w$, then $\mu_{(Y;\Omega)}(w) = 0$.
\endproclaim
\demo{Proof} In this case, we note that the integral (7.15) becomes
$$
\frac{i}{\pi}\int_{\del w} i^*(\omega^j_j) =
\frac{i}{\pi}\int_{w} i^*(d\omega^j_j) =
\frac{i}{\pi}\int_{w} \rho = 0
$$
since $\rho \equiv 0$ for the Calabi-Yau metric.
\qed \enddemo

\proclaim{Corollary 7.6} Suppose that $Y$ is a special
coisotropic submanifold. Then we have
$$
\mu_{(Y;\Omega)}(w)
= \frac{i}{\pi}\int_{\del w} i^*(\omega^a_a).
\tag 7.16
$$
\endproclaim

Note that (7.16) for special Lagrangian submanifolds
reduces to the well-known fact that the Maslov index vanishes 
since in that case $K_Y^{\otimes 2}$ is just a scalar $\C$.

\head 
\n{\bf \S 8. The transverse symplectic curvature in K\"ahler manifolds}
\endhead

In [OP1], we have introduced the notion of $\Pi$-transverse
symplectic curvature $F_\Pi \in \Gamma(\Lambda^2(N^*\FF) \otimes
T\FF)$, which measures non-integrablity of 
the complementary subbundle $G_\Pi$ of  $TY$ 
in a given splitting
$$
TY = T\FF \oplus G_\Pi.
$$
Then the complementary subbunlde $G_\Pi$ is 
integrable, if and only if $F_\Pi = 0$.

In the (almost)-K\"ahler case, we have the canonical
Riemannian splitting
$$
TY = T\FF \oplus (T\FF)^{\perp,J}.
$$
We denote the corresponding curvature as $F = F_Y$. We first recall
the definition

\definition{Definition 8.1 [Definition 4.1, OP1]}
The {\it transverse symplectic curvature} of $Y$ is a $T\FF$-valued
two form on $N_J\FF$ or a section of $\wedge^2(N^*\FF) \otimes T\FF$
defined as follows: Let $\pi: TY \to N_J\FF$ be the orthogonal projection.
For any given $v, \, u \in N_J\FF|_y$, we define
$$
F(v,u): = [X,Y]^\parallel(y)
$$
where $X, \, Y$ be any vector field on $Y$ with $X(y) = v, \, Y(y) = u$
that is normal to the foliation
and $(\cdot)^\parallel$ the component of $(\cdot)$ tangential
to the null-foliation.
\enddefinition

It was proven in [OP1] that this is well-defined as an element
$$
F \in \Gamma(\Lambda^2(N^*\FF) \otimes T\FF).
$$
We will compute the components of the intrinsic curvature
$F$ with respect to those of the second fundamental form of
the embedding $Y \subset (X,g)$ where $g$ is
the associated K\"ahler metric of $X$.

Let $\{e_1, \cdots, e_n, f_1, \cdots, f_n\}$ be an orthonormal frame
of $TX$ adapted to $Y$ and $\omega^i_j = \alpha^i_j + i\beta^i_j$
be the associated Hermitian frame. The first structure equation
is
$$
d\theta^\ell = - \omega^i_j \wedge \theta^j
$$
or equivalently
$$
\cases
de_\ell^*  = - \alpha^\ell_j \wedge e_j^* + \beta_j^\ell \wedge f_j^* \\
df_k^*  = - \beta^k_j \wedge e_j^* - \alpha^k_j \wedge f_j^*.
\endcases
$$
As in section 5, we have
$$
\align
\beta_j^\alpha & = A^\alpha_{j\ell} e_\ell^* + B^\alpha_{jb}f_b^* \\
\alpha_b^\alpha & = C^\alpha_{bj} e_j^* + D^\alpha_{bc}f_c^*
\endalign
$$
where the second fundamental form of $Y \subset (X,g)$ is given by
$S = S^\alpha f_\alpha$ with
$$
S^\alpha = \Big( \matrix A^\alpha_{\ell j}, & C^\alpha_{jb} \\
B^\alpha_{bj}, & D^\alpha_{cb}\endmatrix \Big).
$$

Now we compute $F$. A straightforward calculation leads to
$$
\align
F(e_a, e_b) & = [e_a, e_b]^\parallel = e^*_\alpha([e_a,e_b]) e_\alpha
= - de^*_\alpha(e_a, e_b) e_\alpha \\
& =(-\alpha^\alpha_j(e_a)\delta_b^j + \alpha^\alpha_j(e_b)\delta_a^j)e_\alpha\\
& =(-\alpha^\alpha_b(e_a) + \alpha^\alpha_a(e_b))e_\alpha\\
& = (-C^\alpha_{ba} + C^\alpha_{ab})e_\alpha
\endalign
$$
and similarly
$$
F(e_a, f_b) = (-A^\alpha_{ba} - D^\alpha_{ba})e_\alpha
$$
$$
F(f_a, f_b) = (-B^\alpha_{ba} + B^\alpha_{ab})e_\alpha.
$$
This proves the following formula
$$
\aligned
F & =\Big(\frac{1}{2}(C^\alpha_{ab}
- C^\alpha_{ba})e_\alpha\Big) e^*_a \wedge e^*_b
+\Big(\frac{1}{2}(-B^\alpha_{ba} + B^\alpha_{ab})e_\alpha\Big)
f^*_a \wedge f^*_b \\
& \quad + \Big((-D^\alpha_{ba} - A^\alpha_{ba})e_\alpha\Big) 
e^*_a \wedge f^*_b.
\endaligned
\tag 8.1
$$
In terms of the unitary frame, we can also write
$$
F = F^{2,0} + F^{1,1} + F^{0,2}
\tag 8.2 
$$
where
$$
\align
F^{2,0} & =  \Big(\frac{1}{4}
(C^\alpha_{ab} - C^\alpha_{ba})e_\alpha\Big)\theta^a\wedge\theta^b \tag 8.3\\
F^{1,1} & = \Big(\frac{i}{2}
\Big((-D^\alpha_{ab} - A^\alpha_{ba}) e_\alpha\Big) 
\theta^a\wedge\overline \theta^b \tag 8.4\\
F^{0,2} & = \Big(\frac{1}{4}
(C^\alpha_{ab} - C^\alpha_{ba})e_\alpha \Big)
\overline\theta^a\wedge\overline \theta^b.
\tag 8.5
\endalign
$$
In particular, we have obtained the formula for the {\it symplectic
transverse mean curvature} $\rho^{trans}_Y$ defined in [OP1]
$$
\rho^{trans}_Y = (- A^\alpha_{aa} - D^\alpha_{aa}) e_\alpha.
\tag 8.6
$$

We summarize the above calculation into

\proclaim{Theorem 8.2} \roster
\item
$F$ is of type $(1,1)$ if and only if
$$
C^\alpha_{ab} = C^\alpha_{ba}
\tag 8.7
$$
in addition to (5.7), in which case we have
$$
F = F^{1,1} = \Big(\frac{i}{2} (-D^\alpha_{ab} - A^\alpha_{ba}) e_\alpha\Big)
\theta^a\wedge\overline \theta^b.
$$
\item We have $F = 0$ if and only if the second
fundamental form $S$ satisfies
$$
C^\alpha_{ab} = C^\alpha_{ba}, \quad   A^\alpha_{ba} = -D^\alpha_{ba}
\tag 8.8
$$
in addition to (5.7). 
\endroster
\endproclaim

Now we consider the hypersurface case in detail.

\definition{Example 8.3}
Let $(X,\omega,J)$ be any almost K\"ahler manifold and
consider the hypersurface $Y \subset X$ defined by
$$
Y = \{x \in M \mid \rho(x) = 1 \} \tag 8.9
$$
for a smooth function $\rho: X \to \R$ such that 
$$
|\nabla \rho(x)|_g \equiv 1 \tag 8.10
$$
for any $x \in Y$, where $\nabla \rho$ is the gradient of
$\rho$ with respect to the associated metric $g = \omega(\cdot, J\cdot)$.
We denote by $X_\rho = J\nabla \rho$ is the Hamiltonian vector field
of $\rho$. It follows that $X_\rho$ is tangent to $Y$ and
$$
T\FF = \text{span}_\R\{X_\rho\} \tag 8.10
$$
and $N_J\FF \subset TY$ is nothing but the distribution of the maximal
complex or $J$-invariant subspace of $TY$.

By the definition of $F$, we have
$$
\align
F(X,Y) & = [X,Y]^\parallel = \langle [X,Y], X_\phi\rangle_g X_\rho \\
& = - \omega(J[X,Y],X_\rho)X_\rho = d\rho\circ J ([X,Y]) X_\rho
\tag 8.11
\endalign
$$
for any $X, \, Y \in N_J\FF$.
If we denote
$$
d^c\rho = - d\rho\circ J \tag 8.12
$$
as usual, we have
$$
N_J\FF = \ker d\rho \cap \ker d^c \rho. \tag 8.13
$$
It then follows that
$$
- d^c\rho([X,Y]) = dd^c(X,Y) \tag 8.14
$$
for any $X, \, Y \in N_J\FF$. Combining (8.11) and (8.14), we have 
proved
$$
F(X,Y) = dd^c\rho(X,Y) X_\rho. \tag 8.15
$$
In general, the two forms $dd^c$ on the complex vector
bundle $N_J\FF \otimes \C$ will not be of type $(1,1)$,
unless the almost complex structure is integrable. On the other
hand, in the K\"ahler case, the following idenity is well-known
$$
dd^c\rho = 2i \del \delbar \rho
$$
whose restriction to $N_J\FF$ is the well-known Levi form on the
hypersurface. Therefore in the integrable K\"ahler case,
the transverse symplectic curvature $F$ is automatically
of type (1,1) and the two form
$$
\langle F, X_\phi \rangle |_{N_J\FF}
$$
reduces to the Levi form of the hypersurface.
The case $F = 0$ corresponds to the hypersurface that is Levi-flat
and the associated foliation of $N_J\FF$ in that case is nothing but
the foliation by the maximally complex submanifolds in $Y$.
\enddefinition

Motivated by this example, we now introduce the following definition

\definition{Definition 8.4} We say that a pre-K\"ahler manifold
$(Y,\omega_Y, J_Y)$ is integrable, if the associated transverse
curvature $F$ is of type $(1,1)$, or equivalently the second fundamental
$S = S^\alpha f_\alpha$ in any adapted frame satisfies 
$$
S^\alpha|_{N_J\FF} = \Big( \matrix A^\alpha, C^\alpha \\
C^\alpha, D^\alpha \endmatrix \Big)
$$
where all $A^\alpha, \, C^\alpha$ and $D^\alpha$ are symmetric 
$k \times k$ matrices.
\enddefinition

We will further study geometry of such structures elsewhere
in the future.

\head{\bf Appendix}
\endhead

\n{\bf A.1. Calculation for the almost K\"ahler case}
\smallskip

In this appendix, we generalize the formula in Theorem 4.3 in the
(non-integrable) almost K\"ahler manifold as in the spirit of the
calculation carried out by Schoen and Wolfson [SW] for the
Lagrangian submanifolds in the almost K\"ahler case.

For this, we need to compare the canonical connection used in
the main part of this paper and the Levi-Civita connection of
the metric
$$
g = \omega(\cdot, J\cdot).
$$
The corresponding first structure equation with respect to the
oriented orthonormal frame
$$
\{e_1, \cdots, e_n, f_1, \cdots, f_n \}
$$
can be obtained by decomposing (3.3) into the real and the imaginary
parts: Writing
$$
\align
\omega^i_j & = \alpha^i_j + i \beta^i_j \\
\tau^\ell_{\bar j} & = \gamma_j^\ell + i \delta_j^\ell
\endalign
$$
turns (3.3) into
$$
\align
de_j^* & = - (\alpha^j_\ell + \gamma^j_\ell) \wedge e^*_\ell
+ (\beta^j_\ell - \delta^j_\ell) \wedge f^*_\ell \tag A.1 \\
df_j^* & = - (\beta^j_\ell + \delta^j_\ell) \wedge e^*_\ell
+ (\alpha^j_\ell - \gamma^j_\ell) \wedge f^*_\ell \tag A.2
\endalign
$$
The symmetry properties of $\omega^i_j$ and $\tau^i_{\bar j}$
$$
\overline{\omega}^i_j = - \omega^j_i, \quad \tau^i_{\bar j}
= \tau^j_{\bar i},
$$
coming from the unitarity of $\omega$ and from (3.6), immediately
imply that (A.1-2) defines a torsion free Riemannian connection
of $g$ which is nothing but the unique Levi-Civita connection of $g$.

Now we apply the same kind of analysis as in section 4 noting that
the definitions of $E$-differential and $E$-connections
depend only on the symplectic structure and the given
coisotropic submanifold $Y \subset (X,\omega)$,
but not on the almost complex structure $J$.

In the non-integrable case, there will be torsion terms appearing
in various places. First, (4.4) is replaced by (A.1-2) and (4.5)
by
$$
0 = df^*_\alpha = - (\beta_j^\alpha + \delta_j^\alpha) \wedge e_j^*
- (\alpha^\alpha_j - \gamma^\alpha_j) \wedge f^*_j \tag A.3
$$
on $Y$. Again by Cartan's lemma, (4.6) is replaced by
$$
\aligned
\beta_j^\alpha + \delta_j^\alpha
& = A_{j\ell}^\alpha e_\ell^* + B_{jb}^\alpha f_b^* \\
\alpha_b^\alpha - \gamma_b^\alpha
& = C_{bj}^\alpha e_j^* + D_{bc}^\alpha f_c^*
\endaligned
\tag A.4
$$
By the same reasoning as in section 4 (see [SW] for a similar
calculation for the Lagrangian case), we still obtain the same
formula
$$
\alpha_Y^\parallel = A^\alpha_{\beta\alpha}e^*_\beta \tag A.5
$$
as (4.11).

On the other hand, the formula (4.18) is replaced by
$$
\beta_\gamma^\alpha + \delta_\gamma^\alpha
= A_{\gamma\mu}^\alpha e_\mu^*  \tag A.6
$$
{\it when restricted to } $T\FF$ and so (4.21)
replaced by
$$
i\rho = d\omega^\beta_\beta + d\omega^b_b =
d(\beta_\beta^\beta - \delta_\beta^\beta) + ^Ecurv(K_Y). \tag A.7
$$
Restricting this to $T\FF$ and combining with (A.5-7), we
have obtained the following coisotropic analog to
[Proposition A.1, SW].

\proclaim{Theorem A.1} Let $(X,\omega, J)$ be almost K\"ahler with
the canonical connection. Let $K = (K^i_j)$ be the curvature two form
and $i\rho = K^j_j$ be the Ricci form of the connection.
Then we have the formula
$$
i d_\FF\alpha_Y^\parallel - ^Ecurv(K_Y) = i^*_\FF(\rho) + d_\FF (\tau)
$$
where $\tau \in \Gamma(T^*\FF)$ is a
one form canonically pulled-back via the embedding $i: Y \to X$ from the
torsion form of the canonical connection of $X$.
\endproclaim

\bigskip

\n{\bf A.2. Criterion for the minimality of the null foliation}
\smallskip

In this section, we study the question when the null foliation
becomes a {\it minimal foliation} in the sense that each leaf is a
minimal submanifold of $Y$.

Note that the form
$$
\nu_Y:= e^*_{k+1} \wedge \cdots \wedge e^*_n \in
\Lambda^{n-k}(E^*) = \Omega^{n-k}(\FF)
$$
is independent choice of the frames $\{e_1, \cdots, e_n, f_1,
\cdots, f_k\}$ adapted to $Y$ and so defines a globally
well-defined $n-k$ form on $Y$ which restricts to the volume form
on each leaf.  We recall from [Ru,Su,Ha] that existence of such form
is a necessary and sufficient condition for the manifold $(Y,
g_Y)$ is foliated by {\it minimal} leaves of dimension $n-k$. In
our case, the following theorem shows that the form $\nu_Y$ plays
the role of such form.

\proclaim{Theorem A.2} A leaf $\Sigma$ of the
null foliation is a minimal submanifold if
and only if the from $\nu_Y$ is relatively closed on $\Sigma$: namely,
$$
d\nu_Y(X_1, \cdots, X_{n-k+1}) = 0
\tag A.8
$$
on $\Sigma$ if the first $n-k$ vector fields $X_i$ are tangent to the leaves.
\endproclaim

\demo{Proof}
With respect to the frame $\{e_1, \cdots, d_n, f_1, \cdots, f_k\}$
on $Y$, the first structure equation becomes
$$
\align
de_\ell^* & = - \alpha_j^\ell \wedge e_j^* + \beta_c^\ell \wedge f_c^*
\tag A.9 \\
df_a^* & = -\beta^a_j \wedge e_j^* - \alpha_b^a \wedge f_b^*. \tag A.10
\endalign
$$
On each leaf $\Sigma$ of the null foliation,
provided by $e_a^* = 0 =f_b^*$, we have
$$
\aligned
0 & = de_a^* = - \alpha^a_\alpha \wedge e_\alpha^* \\
0 & = df_a^* = - \beta^a_\alpha \wedge e_\alpha^*. \\
\endaligned
$$
By applying Cartan's lemma and comparing with (4.6), we have
$$
\alpha^a_\alpha = - \alpha^\alpha_a = C^\alpha_{b\beta}e_\beta^*, \,
\beta^a_\beta = \beta_a^\alpha = A^\alpha_{a\gamma}e_\gamma^*.
\tag A.11
$$
Therefore the second fundamental form $S_\Sigma$
of $\Sigma$ in $Y$ is given by the matrix
$$
S_\Sigma = (C^\alpha_{b\beta} e_\alpha^* \otimes e_\beta^*)e_b
+ (A^\alpha_{a\gamma}e_\alpha^* \otimes e_\gamma^*) f_a
\tag A.12
$$
and its mean curvature vector by
$$
\vec H_\Sigma = C^\alpha_{b\alpha}e_b
+ A^\alpha_{a\alpha} f_a.
$$
Hence the leaf $\Sigma$ is minimal in $Y$ if and only if
$$
C^\alpha_{b\alpha} = 0 = A^\alpha_{a\alpha} \tag A.13
$$
for all $a, \, b = 1, \cdots, k$.

Now we compute the differential $d\nu_Y$
$$
\align
d\nu_Y & = d(e_{k+1}^* \wedge \cdots\wedge e_n^*) \\
& = \sum_j (-1)^j e_{k+1}^* \wedge \cdots \wedge de^*_{k+j} \wedge
\cdots \wedge e_n^* \tag A.14
\endalign
$$
The first structure equation (A.9) and (4.6) can be written as
$$
\align
de_\beta^* & = - \alpha^\beta_a \wedge e_a^* - \alpha^\beta_\alpha
\wedge e_\alpha^* \\
& = - (C^\beta_{bj}e_j^* + D^\beta_{bc} f_c^* )\wedge e_b^*
- \alpha^\beta_\alpha e_\alpha^*
- (A^\beta_{b\ell}e_\ell^* + B^\beta_{bc} f_c^*) \wedge f_b^*.
\endalign
$$
It is then straightforward to derive, by substituting this into
(A.14) and using $\alpha_\alpha^\alpha = 0$,
$$
\align
d\nu_Y & = \sum_{k+1 \leq \alpha \leq n}
(-1)^j e_{k+1}^* \wedge \cdots \wedge (-C^\alpha_{b\beta}
e^*_\beta \wedge e_b^*) \wedge \cdots \wedge e_n^* \\
\quad & = \sum_{k+1 \leq \alpha \leq n}
(-1)^j e_{k+1}^* \wedge \cdots \wedge (A^\alpha_{b\gamma}
e^*_\gamma\wedge f_b^*) \wedge \cdots \wedge e_n^* \\
& = (C^\alpha_{b\alpha}e_b^* - A^\alpha_{b\alpha}f_b^*)\wedge
\nu_Y \tag A.15
\endalign
$$
for the $(n-k +1)$-tuple  $(X_1, \cdots, X_{n-k+1})$ with the
first $n-k$ vector fields $X_i$ are tangent to the leaves. This
finishes the proof. \qed\enddemo

\bigskip

\head{References}\endhead 

\enddocument